\documentclass[12pt]{amsart}
\usepackage{etex}
\usepackage{amssymb}
\usepackage{amsfonts}
\usepackage{amsmath}
\usepackage{array}
\usepackage{rotating}
\usepackage{pstricks}
\usepackage[utf8]{inputenc}

\usepackage{graphicx}
\usepackage{placeins}
\usepackage{pspicture}

\usepackage{tikz}
\usepackage{calrsfs}

\usepackage{skak} 

\usepackage[matrix,ps,xdvi]{xypic} 
\newdimen\vcadre\vcadre=0.1cm 
\newdimen\hcadre\hcadre=0.1cm 
\def\GrTeXBox#1{\vbox{\vskip\vcadre\hbox{\hskip\hcadre%
      $#1$%
   \hskip\hcadre}\vskip\vcadre}}
\def\arx#1[#2]{\ifcase#1 \relax \or%
  \ar @{-}[#2]  \or%
  \ar @2{-}[#2] \or%
  \ar @{--}[#2] \or%
  \ar @2{.}[#2] \or%
  \ar @{~}[#2]  \fi}

\def\arbgx#1#2{
\newdimen\vcadre\vcadre=0.01cm 
\newdimen\hcadre\hcadre=0.01cm 
\xymatrix@R=0.1cm@C=1mm{
 & {\GrTeXBox{\bullet}}\arx1[dl]\arx1[dr]\\
 {\GrTeXBox{#1}} & *{} & {\GrTeXBox{#2}} \\
}
}

\def\arbgxb#1#2#3{
\newdimen\vcadre\vcadre=0.01cm 
\newdimen\hcadre\hcadre=0.01cm 
\xymatrix@R=0.1cm@C=2mm{
 & {\GrTeXBox{\bullet}}\arx1[dl]\arx1[dr]\\
 {\GrTeXBox{#1}} & *{} & {\GrTeXBox{\bullet}}\arx1[dl]\arx1[dr]\\
 & {\GrTeXBox{#2}} & *{} & {\GrTeXBox{#3}} & *{} \\
}
}

\def\arbgxc#1#2#3{
\newdimen\vcadre\vcadre=0.01cm 
\newdimen\hcadre\hcadre=0.01cm 
\xymatrix@R=0.1cm@C=2mm{
 && {\GrTeXBox{\bullet}}\arx1[dl]\arx1[dr]\\
 & {\GrTeXBox{\bullet}}\arx1[dl]\arx1[dr]  & *{} & {\GrTeXBox{#3}} \\
 {\GrTeXBox{#1}} & *{} & {\GrTeXBox{#2}} & *{} \\
}
}

\newtheorem{example}{Example}[section]

\newtheorem{theorem}[example]{Theorem}

\newtheorem{proposition}[example]{Proposition}

\newcommand{\tmmathbf}[1]{\ensuremath{\boldsymbol{#1}}}

\newcommand{\sqx}{\sympawn}

\def\Proof{\noindent \it Proof -- \rm}

\def\qed{\hspace{3.5mm} \hfill \vbox{\hrule height 3pt depth 2 pt width 2mm}
\bigskip}

\def\xx{{\bf x}}
\def\zz{{\bf z}}

\def\std{{\rm std}}

\def\maj{{\rm maj}}

\def\FQSym{{\bf FQSym}}
\def\WQSym{{\bf WQSym}}

\def\Cc{{\sf C}}
\def\DD{{\sf D}}

\def\Lie{{\rm Lie}}
\def\<{\langle}
\def\>{\rangle}

\def\R{{\mathbb R}}

\def\CC{{\mathfrak C}}
\def\C{\operatorname{\mathbb C}}

\def\KK{\operatorname{\mathbb K}}

\def\Hh{{\mathcal H}}
\def\F{{\bf F}}
\def\S{{\bf S}}
\def\G{{\bf G}}

\def\M{{\bf M}}
\def\P{{\bf P}}

\def\SG{{\mathfrak S}}

\def\AA{{\mathcal A}}

\def\Sym{{\bf Sym}}

\def\Ff{{\mathcal F}}

\def\Res{\operatorname{Res}}

\def\c{{\sf c}}

\def\PBT{{\bf PBT}}
\def\shuff#1#2{\mathbin{
\hbox{\vbox{ \hbox{\vrule \hskip#2 \vrule height#1 width 0pt
}%
\hrule}%
\vbox{ \hbox{\vrule \hskip#2 \vrule height#1 width 0pt
\vrule }%
\hrule}%
}}}

\def\shuf{{\mathchoice{\shuff{7pt}{3.5pt}}%
{\shuff{6pt}{3pt}}%
{\shuff{4pt}{2pt}}%
{\shuff{3pt}{1.5pt}}}}%
\def\shuffle{\,\shuf\,}

\def\ZZ{{\mathbb Z}}    
\def\KK{{\mathbb K}\, } 

\def\x{{\sf x}}

\newdimen\Squaresize \Squaresize=14pt
\newdimen\Thickness \Thickness=0.5pt

\def\Square#1{\hbox{\vrule width \Thickness
   \vbox to \Squaresize{\hrule height \Thickness\vss
      \hbox to \Squaresize{\hss#1\hss}
   \vss\hrule height\Thickness}
\unskip\vrule width \Thickness}
\kern-\Thickness}

\def\Vsquare#1{\vbox{\Square{$#1$}}\kern-\Thickness}

\def\qbin#1#2{\begin{bmatrix} #1 \\ #2\end{bmatrix}}
\def\ssbin#1#2{ \big( {\scriptstyle #1 \atop \scriptstyle #2 } \big)}

\def\Tabvrule{\vrule width-0.4pt}       
\def\Tabhrule{\hrule \hrule height-0.4pt} 
\def\Tabstrut{\vrule height2.2ex 
                     depth0.8ex  
                     width0ex    
\relax}

\def\PasCase#1{\omit%
            $\vcenter{\hbox {\vbox to 0.4pt{}}
               \hbox{\makebox[3ex]{\Tabstrut$#1$}}}%
               \Tabvrule$}
\def\PasCasePoint{\PasCase{\cdot}}
\def\DessinCarre#1{%
    \vcenter{\hbox{}\hrule
             \hbox{\vrule\makebox[3ex]{\Tabstrut$#1$}\vrule}\Tabhrule}%
             \Tabvrule}
\def\GenRuban#1{\vcenter{\halign{&$\DessinCarre{##}$\cr#1}}\egroup}

\def\sTabvrule{\vrule width-0.4pt}
\def\sTabhrule{\hrule \hrule height-0.4pt}
\def\sTabstrut{\vrule height1.6ex depth0.6ex width0ex \relax}
\def\sDessinCarre#1{%
    \vcenter{\hbox{}\hrule
             \hbox{\vrule\makebox[2.3ex]%
                  {\sTabstrut$\scriptstyle#1$}\vrule}\sTabhrule}%
             \sTabvrule}
\def\sGenRuban#1{\vcenter{\halign{&$\sDessinCarre{##}$\cr#1}}\egroup}

\def\ruban{%
  \bgroup
  \let\ =\omit
  \let\\=\cr
  \let\x=\times
  \let\.=\PasCasePoint
  \offinterlineskip
  \GenRuban}

\def\sruban{%
  \bgroup
  \let\ =\omit
  \let\x=\times
  \let\\=\cr
  \offinterlineskip
  \sGenRuban}

\def\arbuga{\begin{picture}(2,4)\cerg{1}3\cerp{1}1
 \put(1,3){\Line( 0,-2)}
\end{picture}}

\def\arbdga{\begin{picture}(2,6)\cerg15\cerp13\cerp11
 \put(1,5){\Line( 0,-2)}
 \put(1,3){\Line( 0,-2)}
\end{picture}}

\def\arbdgb{\begin{picture}(2,4)\cerg15\cerp13\cerp11
 \put(1,5){\Line( 0,-2)}
 \put(1,3){\Line( 0,-2)}
\end{picture}}

\def\arbdgb{\begin{picture}(3,4)\cerg23\cerp11\cerp31
 \put(2,3){\Line( 1,-2)}
 \put(2,3){\Line(-1,-2)}
\end{picture}}

\newdimen\vcadre\vcadre=0.2cm 
\newdimen\hcadre\hcadre=0.2cm 

\def\cerp#1#2{\put(#1,#2){\circle*{0.7}}}
\def\cerg#1#2{\put(#1,#2){\circle*{1}}}

\def\arbtga{\begin{picture}(3,4)\cerg23\cerp{.5}1\cerp21\cerp{3.5}1
 \put(2,3){\Line(-1.5,-2)}
 \put(2,3){\Line( 0,-2)}
 \put(2,3){\Line( 1.5,-2)}
\end{picture}}

\def\arbtgb{\begin{picture}(3,6)\cerg25\cerp13\cerp33\cerp11
 \put(2,5){\Line(-1,-2)}
 \put(2,5){\Line( 1,-2)}
 \put(1,3){\Line( 0,-2)}
\end{picture}}

\def\arbtgd{\begin{picture}(3,6)\cerg25\cerp13\cerp33\cerp31
 \put(2,5){\Line(-1,-2)}
 \put(2,5){\Line( 1,-2)}
 \put(3,3){\Line( 0,-2)}
\end{picture}}

\def\arbtgc{\begin{picture}(3,6)\cerg25\cerp23\cerp11\cerp31
 \put(2,5){\Line( 0,-2)}
 \put(2,3){\Line( 1,-2)}
 \put(2,3){\Line(-1,-2)}
\end{picture}}

\def\arbtge{\begin{picture}(3,8)\cerg27\cerp25\cerp23\cerp21
 \put(2,7){\Line( 0,-2)}
 \put(2,5){\Line( 0,-2)}
 \put(2,3){\Line( 0,-2)}
\end{picture}}

\def\cerp#1#2{\put(#1,#2){\circle*{0.7}}}
\def\cerg#1#2{\put(#1,#2){\circle*{1}}}

\title[]{Tree expansions of some Lie idempotents}

\author[] {Frédéric Menous, Jean-Christophe Novelli and Jean-Yves Thibon}

\date{\today}

\begin{document}

\begin{abstract}
	We prove that the Catalan Lie idempotent $D_n(a,b)$, introduced in
	[Menous {\it et al.}, Adv.  Appl. Math. 51 (2013), 177] can be 
	refined by introducing $n$ independent parameters $a_0,\ldots,a_{n-1}$ 
	and that the coefficient of each monomial is itself a Lie idempotent
	in the descent algebra. These new idempotents are multiplicity-free
	sums of subsets of the Poincaré-Birkhoff-Witt basis of the Lie module.
	These results are obtained by embedding noncommutative symmetric functions into the
	dual noncommutative Connes-Kreimer algebra, which also allows us to interpret, and rederive in a simpler way,
	Chapoton's results on a two-parameter tree expanded series.
\end{abstract}

\maketitle

\small
\tableofcontents
\normalsize
\newpage

\section{Introduction}

Lie idempotents are idempotents of the symmetric group algebra which act on words as projectors
onto the free Lie algebra. Thus, they are in particular elements of the Lie module $\Lie(n)$,
spanned by complete bracketing of standard words, such as $[[1,3],[2,4]]$,  
which can be represented as complete binary trees with leaves labelled $1,2,\ldots,n$.

Of course, these elements are not linearly independent, but the trees  such that for each
internal node, the smallest label is in the left subtree, and the greatest
label is in the right subtree do form a basis, called the Poincaré-Birkhoff-Witt basis \cite{ST}.
Such labellings are said to be admissible. These basis  elements are denoted by $t(\sigma)$, where $t$
is a complete binary tree, and $\sigma$ the permutation obtained by reading its leaves from left to right.

The direct sum $\Lie = \bigoplus_{n\ge 0}\Lie(n)$ can be interpreted as the operad ${\mathcal Lie}$.
It is also a Lie algebra for the Malvenuto-Reutenauer convolution product of permutations,
which allows us to regard it as contained into $\FQSym$, a permutation $\sigma$ being intepreted as the basis
element $\G_\sigma$. Then, it is (strictly) contained in the primitive Lie algebra of $\FQSym$.

It turns out that the elements $\c_t$, defined for complete binary trees $t$ by
the sum over  admissible labellings
\begin{equation}
\c_t = \sum_{\sigma\ {\rm admissible}}t(\sigma)
\end{equation}
span a Lie subalgebra $\CC$  of $\Lie$, which might be called the Catalan Lie algebra.

It has proved convenient to label its basis elements by plane trees instead of binary trees:
we set $C_T=\c_t$ where the plane tree $T$ is the  right-branch rotation of the
incomplete binary tree $t'$ obtained by removing the leaves of $t$ (so that the maximal element
of the Tamari order is the corolla).

This provides us with elements $C_T$ of $\FQSym$, labelled by plane trees. 
The noncommutative Connes-Kreimer algebra $\Hh_{NCK}$ \cite{Foi00,Foi01} is the free associative algebra generated
by variables $Y_T$ indexed by plane trees, endowed with the coproduct of admissible cuts.
Its basis $Y_F$ is therefore indexed by plane forests $F$.
Its dual $\Hh_{NCK}^*$ admits a natural embedding into $\FQSym$, and if $X_F$ denotes the dual basis
of $Y_F$, it turns out that
\begin{equation}
	C_T = \sum_{T'\le T}X_{T'}
\end{equation}
where the sum is over the Tamari order. Moreover,
if one denotes by $\tau=\bar T$ the underlying non-plane rooted tree of $T$
the sums
\begin{equation}
x_\tau = |{\rm Aut}(\tau)|\sum_{\bar T=\tau}X_T
\end{equation}
span a sub-preLie algebra, which is free on the generator $x_\bullet$, and 
$x_\tau$ coincides with its Chapoton-Livernet basis.

The aim of this paper is to investigate the expansions in the $X$ and $C$ bases of various
noncommutative symmetric functions, regarded as elements of $\FQSym$.
Our first result concerns the family of Catalan idempotents $D_n(a,b)$.
Originally introduced as noncommutative symmetric functions on the ribbon basis in \cite{MNT},
these elements were identified in \cite{FMNT} as simple weighted sums of the basis $C_T$.
However, the calculations of \cite{FMNT} are rather tricky, and it is by no means obvious
that such sums belong to the descent algebra. We present here a new approach, relying on the
Birkhoff factorisation of a simple character of $QSym$ with values in an algebra of Laurent series.
This approach produces immediately the expansion on $X_F$ of a grouplike series $\sigma_{a(z)}^+$ which by definition
belongs to the descent algebra. It is then relatively straightforward to check that
the original  Catalan idempotents are obtained by choosing $a(z)=\frac{a}{z}+\frac{b}{1-z}$ and taking the residue,
the general case giving rise to new refined idempotents indexed by partitions of $n-1$.

Finally, we show how the embedding of $\Hh_{NCK}^*$ into $\FQSym$ can be used to determine
the $X$-expansion of various noncommutative symmetric functions, including the Eulerian idempotents and
the two-parameter series of Chapoton \cite{Cha3}

This paper is a continuation of \cite{FMNT}, to which the reader is referred for background and notation.

{\bf Acknowlegements. } This research has been partially supported by the project CARPLO
of the Agence Nationale de la recherche (ANR-20-CE40-0007).

\section{The noncommutative Connes-Kreimer Hopf algebra $\Hh_{NCK}$}

The noncommutative Connes-Kreimer Hopf algebra $\Hh_{NCK}$, introduced by Foissy \cite{Foi00,Foi01},
is as a graded vector space spanned by plane forests, the degree being the number of nodes.
We denote by $Y_F$ its natural basis indexed by forests:
\begin{equation}
\Ff=\lbrace \emptyset,\bullet, \arbuga,\bullet \bullet, \arbdga ,\arbdgb ,\bullet \arbuga,\arbuga \bullet, \bullet \bullet \bullet,\dots \rbrace
\end{equation}
 It is then freely generated by variables $Y_T$ indexed by plane trees. 
The product is concatenation, and the coproduct is
given by admissible cuts, which can be conveniently defined directly for the iterated coproducts in terms
of labellings. 

Trees will be drawn with the root at the top in this paper. The canonical labelling of a tree
is obtained by visiting it in postorder, so that the labels of each subtree form an interval,
with the maximum at its root. 

\medskip
{\footnotesize
For instance, for the tree $\arbtgb$, we get the labelling
{ \newcommand{\nodea}{\node[draw,circle] (a) {$2$}
;}\newcommand{\nodeb}{\node[draw,circle] (b) {$1$}
;}\newcommand{\nodec}{\node[draw,circle] (c) {$3$}
;}\newcommand{\noded}{\node[draw,circle] (d) {$3$}
;}\newcommand{\nodee}{\node[draw,circle] (e) {$6$}
;}\newcommand{\nodef}{\node[draw,circle] (f) {$4$}
;}\newcommand{\nodeg}{\node[draw,circle] (g) {$5$}
;}
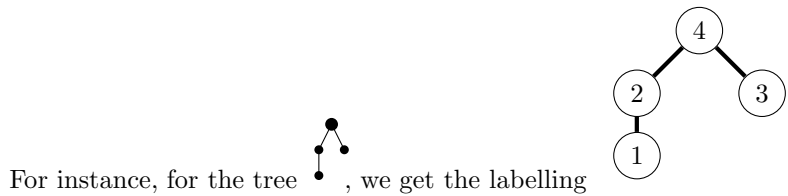
\begin{tikzpicture}[auto]
\matrix[column sep=.2cm, row sep=.2cm,ampersand replacement=\&]{
         \& \nodef  \&         \\ 
 \nodea  \& 	 \&  \nodec      \\ 
	\nodeb  \&         \&  \\
};

\path[ultra thick, black] (f) edge (a)
	(f) edge (c) 
	(a) edge (b);
\end{tikzpicture}}
}

\medskip
A forest $F$ is similarly labelled, by first grafting it on a common root -- that is considering the tree $T=B_+(F)$ -- labelling $T$ 
and removing this labelled root afterwards.

Such a labelled forest is regarded as the Hasse diagram of a poset.

With this labelling the $r$-iteraded coproduct of an element $Y_F$ of degree n can be described as follows:
\begin{equation}
	\Delta^r Y_F = \sum_{u\in C(F)\cap [r]^n}Y_{F_{(1)}}\otimes Y_{F_{(2)}}\cdots\otimes Y_{F_{(r)} }
\end{equation}
where $C(F)$ is the set of words such that
$i<_Fj\Rightarrow u_i\le u_j$,
and  $F_{(i)}$ is the  restriction of $F$ to vertices labelled $i$.

\medskip
{\footnotesize
For instance, for the previous tree $T=\arbtgb$ and $r=2$, 
\begin{equation}
C(T)\cap \lbrace 1,2 \rbrace^4=\lbrace 2222, 2212, 1222, 1122, 1212, 1112, 1111  \rbrace
\end{equation}
give the coproduct:
\begin{equation}
\Delta\arbtgb = 1\otimes\arbtgb + 
\bullet\otimes \arbdga +
\bullet\otimes \arbdgb 
	+\arbuga\otimes\arbuga
	+\bullet\bullet\otimes \arbuga+
	\arbuga\bullet\otimes\bullet+
\arbtgb\otimes1 
\end{equation}
}

\medskip
As it will be useful in the following sections, let us also recall here  the polish code of a plane forest is obtained by labelling each node by the number of its descendants, and traversing it in prefix order. For the previous  tree $T$ we get the polish code $2100$ and also its reverse polish code $0012$.

It has been shown in \cite[3.5]{FNT} that
$\Hh_{NCK}$ admits an embedding $\pi$ into $\WQSym$. 
It is actually an embedding into $\FQSym$, given  by  $F\mapsto \Gamma_F(A)$, where $\Gamma_P(A)$ denotes
the free generating function of a poset \cite{NCSF6}, that is, the sum of its linear extensions
\begin{equation}
	\Gamma_P(A) = \sum_{\sigma\in L(P)}\F_\sigma \in\FQSym = \sum_{u\in C(F)}\M_u,
\end{equation}
where $C(F)$ is the set of packed words $u$ such that $i<_F j\Rightarrow u_i\le u_j$.
Indeed, the linear extensions of a poset are precisely those permutations $\sigma$
such that $i<_P j\Rightarrow \sigma^{-1}(i)<\sigma^{-1}(j)$.

The linear extensions of such a labelled  forest  form an initial interval of the right weak order \cite{BW2}.

\bigskip
{\footnotesize
For example,
\begin{equation}
	\Gamma_{\arbtgb} = \F_{3124}+\F_{1324}+\F_{1234}=\S^{2314}=\check\S^{3124}
\end{equation}
where \cite[(6.4), (6.12)]{NCSF7}
\begin{equation}
	\S^\sigma=\sum_{\tau\le_L \sigma}\G_{\tau}=:\check \S^{\sigma^{-1}}.
\end{equation}
Then, $Y_F$ can be identified with $\Gamma_F=\check\S^{\sigma_F}$, where $\sigma_F$
is the maximal linear extension of $F$.
For example
\begin{equation}
\Gamma_{\arbuga}\Gamma_\bullet = \S^{12}\S^1 = \S^{231} = \check\S^{312} = \Gamma_{\arbuga\bullet}.
\end{equation}
As for the coproduct, $\Gamma_{\arbtgb}=\check\S^{3124}$, and 
{\small
\begin{equation}
\Delta\check\S^{3124}=
	1 \otimes \check\S^{3124} + \check\S^{1} \otimes \check\S^{123} +  \check\S^{1} \otimes \check\S^{213} +
	\check\S^{1 2} \otimes \check\S^{1 2} + \check\S^{21} \otimes \check\S^{12} + \check\S^{312} \otimes \check\S^{1} + 
 \check\S^{3124} \otimes 1,
\end{equation}
}
which corresponds term by term to
\begin{equation}
\Delta\arbtgb = 1\otimes\arbtgb + 
\bullet\otimes \arbdga +
\bullet\otimes \arbdgb 
	+\arbuga\otimes\arbuga
	+\bullet\bullet\otimes \arbuga+
	\arbuga\bullet\otimes\bullet+
\arbtgb\otimes1 
\end{equation}
}

Indeed, the coproduct formula \cite[(6.13)]{NCSF7}
\begin{equation}
\Delta\check\S^\sigma=
\sum_{u\cdot v\le\sigma}\<\sigma|u\shuffle v\>
\check\S^{\std(u)}\otimes\check\S^{\std(v)}\,,
\end{equation}
(sum over pairs of complementary subwords whose concatenation is smaller than $\sigma$
in the right weak order) 
implies that if a value $\sigma_i$ goes into $v$, all greater values on its right must
also go into $v$, so as not to create new inversions. Thus, the word $u$ and $v$ encode
admissible cuts.

\section{Dual noncommutative Connes-Kreimer algebra $\Hh_{NCK}^*$}

\subsection{An embedding of $\Sym$, and its dual}
Let $X_F$ be the dual basis of $Y_F$. According to our description of the coproduct of $Y_F$,
the coefficient of $X_F$ in the product $X_{F'}X_{F''}$ is equal to the number of labellings of $F$
by words over $\{1,2\}$, nondecreasing from bottom to top, and  such that
 $F_{(1)}=F'$ and $F_{(2)}=F''$.

The coproduct of $X_F$ is deconcatenation, so that trees $X_T$ are primitive.
The elements
\begin{equation}\label{eq:SymInH}
\Lambda_n := X_{\bullet\bullet\cdots\bullet} \quad\text{($n$ vertices) and }\quad S_n:=\sum_{|F|=n}X_F
\end{equation}
form sequences of divided powers, and both define the same embedding of $\Sym$ into
$\Hh_{NCK}^*$.
One easily checks that, indeed,
\begin{equation}
\left(\sum_{n\ge 0}(-1)^n X_{\bullet\bullet\cdots\bullet}\right)^{-1}=\sum_{n\ge 0}\sum_{|F|=n}X_F.
\end{equation}
\bigskip
{\footnotesize
Representing trees by their Polish codes, we have:
\begin{align*}
R_{1 1} & = X_{00}\\
R_{2} & = X_{00}+ X_{10}\\
R_3 &= X_{000}+X_{100}+X_{010}+X_{200}+X_{110}\\
R_{21} &= 2X_{000}+ X_{100}+ X_{010}\\
R_{12}&=   2X_{000}+ X_{100}+ X_{010}+X_{200}\\
R_{111}&=X_{000}\\
R_4&= X_{0000} + X_{0010} + X_{0100} + X_{1000} + X_{1010} + X_{0200} + X_{2000}\\
& + X_{1100} + X_{0110} + X_{1110} + X_{1200} + X_{2100} + X_{2010} +  X_{3000}  \\
R_{31}&= (X_{1100}+ X_{0110})+ 2(X_{1000}+ X_{0100}+ X_{0010})\\
&\ + (X_{2000}+ X_{0200})+ X_{1010}+ 3X_{0000}\\
R_{22}&= (X_{1100}+ X_{0110})+ 3(X_{1000}+ X_{0100}+ X_{0010})\\
&\ + 2(X_{2000}+ X_{0200})+ 2X_{1010} + 2X_{3000}+(X_{2100}+X_{2010})+5X_{0000}\\
R_{13}&= (X_{1100}+ X_{0110})+ 2(X_{1000}+ X_{0100}+ X_{0010})\\
&\ + 2(X_{2000}+ X_{0200})+ X_{1010} + 2X_{3000}+(X_{2100}+X_{2010})+X_{1200}+3X_{0000}\\
R_{211}&= (X_{1000}+ X_{0100}+ X_{0010})+3X_{0000}\\
R_{121}&=  2(X_{1000}+ X_{0100}+ X_{0010})+ (X_{2000}+ X_{0200})+ X_{1010}+5X_{0000}\\
R_{112}&=(X_{1000}+ X_{0100}+ X_{0010})+(X_{2000}+ X_{0200})+X_{3000}+3X_{0000}\\
R_{1111}&= X_{0000}
\end{align*}
}

\begin{proposition}
The coefficient of $ X_F$ in  $R_I$ is equal to the number of linear extensions of $F$
which are of ribbon shape $I$.
\end{proposition}

\Proof 
The product rule of the $X$-basis implies immediately that the coefficient 
$\< Y_F,S^I\>$ of $X_F$ in $S^I$
is the number of nondecreasing labellings of $F$ by words of evaluation $I$,
which is also the coefficient of $M_I$ in $\Gamma_F(X)$.
The dual map of the embedding of $\Sym$ into $\Hh_{NCK}^*$ is an epimorphism
$\pi:\ \Hh_{NCK}\rightarrow QSym$, and  $\< Y_F,S^I\>$ is also the coefficient
of $M_I$ in $\pi(Y_F)$. 

Thus, $\pi$ is the restriction of the canonical projection
(commutative image) of $\FQSym$ onto $QSym$, so that the coefficient of $X_F$ in $R_I$
is equal to the coefficient of $F_I$ in $\pi(\Gamma_F)$, which is by definition the
number of linear extensions of $F$ of ribbon shape $I$.
\qed

\medskip
{\footnotesize
For instance, if $T=\arbtgb=2100$, we have $\Gamma_{T} = \F_{3124}+\F_{1324}+\F_{1234}$ and $(3124)$, $(1324)$, $(1234)$ have respective ribbon shape $13$, $22$ and $4$,
so that $X_{2100}$ appears with a coefficient $1$ in $R_{13}$, $R_{22}$ and $R_{4}$.
}
\medskip

The product rule of the $X$-basis also implies that the coefficient of $X_F$ in $\Lambda^I$ is the number of strictly
decreasing labellings of $F$ by words of evaluation $I$.

We have therefore proved:
\begin{theorem}\label{th:coeffX}
The coefficient of $X_F$ is $\sigma_1(XA)$ is
\begin{equation}
\sum_I\<Y_F,S^I\>M_I(X) = \Gamma_F(X)
\end{equation}
and that of $X_F$ in $\lambda_1(XA)$ is
	\begin{equation}\label{eq:chi}
\sum_I\<Y_F,\Lambda^I\>M_I(X) = (-1)^{|F|}\Gamma_F(-X) =: \chi_F(X).
\end{equation}
\end{theorem}
\qed

Recall from \cite{NCSF2} that the involution $X\mapsto -X$ of $QSym$ is  the adjoint of $A\mapsto -A$ in $\Sym$,
so that 
\begin{equation}\label{eq:def-X}
M_I(-X)= (-1)^{\ell(I)}\sum_{J\le I}M_J(X)\quad\text{et}\quad F_I(-X)=(-1)^{|I|}F_{\bar I^\sim}(X).
\end{equation}

\subsection{Dendriform structure}
One can also describe the product $X_{F_1}X_{F_2}$  in the following way.
Endow $F_1$ with its canonical labelling, and $F_2$ with its canonical labelling shifted by
the number $n_1$ of vertices of $F_1$. Then, the coefficient of $X_F$ in the product is equal
to the number of standard decreasing (from the roots towards the leaves) labellings of $F$ whose
restriction to $[1,n_1]$ is $F_1$, and whose restriction to 
$[n_1+1,n_1+n_2]$ is $F_2$. 

This allows to define dendriform half-products:
$X_{F_1}\prec X_{F_2}$ consists of the forests whose first label in the postfix order is $n_1+1$,
and $X_{F_1}\succ X_{F_2}$ of those whose first label is $1$.
In particular,
\begin{equation}
X_T\prec X_F = X_{FT},
\end{equation}
Actually, $\Hh_{NCK}^*$ can be identified with the Loday-Ronco Hopf algebra $\PBT$ \cite{Foi00,Foi01}. 
It is easily seen that its dendriform product are induced by those of $\FQSym$, so that 
the coefficient of 
$X_F$ in $\P_t$ is equal to the number of linear extensions of $F$ whose decreasing
tree has shape $t$. 

\subsection{PreLie and brace structures}
The product rule shows that for two trees, 
$X_{T_1T_2}$ et $X_{T_2T_1}$ have the same coefficient in
$X_{T_1}X_{T_2}$. 
Thus, $[X_{T_1},X_{T_2}]$ is a linear combination of trees, and the primitive Lie algebra 
admits the $X_T$ as basis.

The commutator $[X_{T_1},X_{T_2}]$ is the difference between the sum of the  $X_T$ obtained by grafting
$T_1$ on a node of $T_2$ and that of those obtained by grafting $T_2$ on a node of  $T_1$.
If one denotes by $X_{T_1}\triangleright X_{T_2}$ the first sum, $\triangleright$ defines then  a right preLie product.
There is also a left preLie product. $a\triangleleft b=b\triangleright a$.

\bigskip
{\footnotesize
For example,
\begin{equation}
	[X_\bullet,X_{\arbuga} ] = 2X_{\arbdgb} = X_\bullet\triangleright X_{\arbuga}-X_{\arbuga}\triangleright X_\bullet.
\end{equation}
}

The preLie algebra generated by $X_\bullet$ is free. For a non-plane rooted tree $\tau$, we set
\begin{equation}
x_\tau = |{\rm Aut}(\tau)|\sum_{\bar T= \tau}X_T
\end{equation}
(where $\bar T$ means that the rooted tree obtained by forgetting the planar structure of $T$ is $\tau$)
gets identified with the Chapoton-Livernet basis of the free preLie algebra on one generator.

The brace product which extends the preLie product and induces the associative product is \cite{Foi3}
\begin{equation}
	\< X_{T_1\cdots T_r}, X_T\>_\triangleright = \sum_{T'}X_{T'}
\end{equation}
where the sum runs over all trees $T'$ obtained by grafting
$T_1,\ldots,T_r$ on nodes of $T$, respecting their order.
One has then,
writing $B(X_F)$ for $X_{B_+(F)}$,
$B(X_FX_{F'})=\<X_F, B(X_{F'})\>$
and
\begin{equation}
	\<X_F, \<X_{F'},X_T\>\> =\<X_FX_{F'},X_T\>.
\end{equation}
The primitive Lie algebra of 
$\Hh_{NCK}^*$ is the free brace algebra on one generator.

In terms of the dendriform operations, the preLie product is  $x\triangleright y =a\succ b-b\prec a$, and
one has then as usual $\Lambda_n=X_\bullet\prec\Lambda_{n-1}$
and $S_n=S_{n-1}\succ X_\bullet$.

\subsection{A quotient of $\FQSym$}
Let $\M_\sigma$ be the dual basis of $\S^\sigma$.
The above embedding of  
$\Hh_{NCK}$ into $\FQSym$ allows to identify $X_F$
with the image of $\M_{\sigma_F^{-1}}$ in the quotient of  $\FQSym$ by the
relations $\M_\sigma\equiv 0$ if $\sigma$ contains the pattern  $132$.

\bigskip
{\footnotesize
For example,
\begin{equation}
X_{\arbuga}X_{\arbuga} = 2X_{\arbuga\arbuga}+X_{\arbtgb}+X_{\arbtgd}+X_{\arbtge}
\end{equation}
and
\begin{equation}
\begin{aligned}
\M_{12}\M_{12}& =
	\M_{1 2 3 4} + 2\M_{{\bf 1 3 2} 4} + \M_{{\bf 1 3} 4 \bf 2} + \M_{{\bf 1 4 2} 3}\\ 
&+ \M_{2 3 1 4} 
	+ \M_{{\bf  2 4} 1 \bf 3} + \M_{ 3 1 2 4} + \M_{ 3 \bf 1 4 2} + 2\M_{3 4 1 2}\\
&\equiv   2\M_{3 4 1 2} +  \M_{3124} +  \M_{2314} +\M_{1 2 3 4}\\
& = 2\M_{\sigma_{\arbuga\arbuga}^{-1}}+\M_{\sigma_{\arbtgb}^{-1}}+\M_{\sigma_{\arbtgd}^{-1}}+\M_{\sigma_{\arbtge}^{-1}}
\end{aligned}
\end{equation}
}

To reconstruct the forest $F$ from its maximal linear extension $\sigma_F$, one must construct
the binary search tree of its mirror image $\overline{\sigma_F}$ and take its right branch.rotation 


\bigskip
{\footnotesize
For example, the tree $3100200$ 

{ \newcommand{\nodea}{\node[draw,circle] (a) {$7$}
;}\newcommand{\nodeb}{\node[draw,circle] (b) {$2$}
;}\newcommand{\nodec}{\node[draw,circle] (c) {$1$}
;}\newcommand{\noded}{\node[draw,circle] (d) {$3$}
;}\newcommand{\nodee}{\node[draw,circle] (e) {$6$}
;}\newcommand{\nodef}{\node[draw,circle] (f) {$4$}
;}\newcommand{\nodeg}{\node[draw,circle] (g) {$5$}
;}
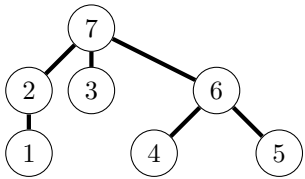
\begin{tikzpicture}[auto]
\matrix[column sep=.2cm, row sep=.2cm,ampersand replacement=\&]{
         \& \nodea  \&         \&         \&         \\ 
 \nodeb  \& \noded  \&         \& \nodee  \&         \\ 
 \nodec  \&         \& \nodef  \&         \& \nodeg  \\
};

\path[ultra thick, black] (b) edge (c)
	(e) edge (f) edge (g)
	(a) edge (b) edge (d) edge (e);
\end{tikzpicture}}

has $\sigma_T = 5463127$, and the binary search tree of  $\overline{\sigma_T }$ is

{ \newcommand{\nodea}{\node[draw,circle] (a) {$7$}
;}\newcommand{\nodeb}{\node[draw,circle] (b) {$2$}
;}\newcommand{\nodec}{\node[draw,circle] (c) {$1$}
;}\newcommand{\noded}{\node[draw,circle] (d) {$3$}
;}\newcommand{\nodee}{\node[draw,circle] (e) {$6$}
;}\newcommand{\nodef}{\node[draw,circle] (f) {$4$}
;}\newcommand{\nodeg}{\node[draw,circle] (g) {$5$}
;}
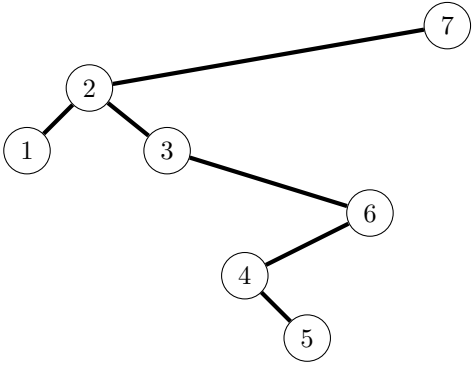
\begin{tikzpicture}[auto]
\matrix[column sep=.2cm, row sep=.2cm,ampersand replacement=\&]{
         \&         \&         \&         \&         \&         \&         \&         \&         \& \nodea  \&         \\ 
         \& \nodeb  \&         \&         \&         \&         \&         \&         \&         \&         \&         \\ 
 \nodec  \&         \&         \& \noded  \&         \&         \&         \&         \&         \&         \&         \\ 
         \&         \&         \&         \&         \&         \&         \& \nodee  \&         \&         \&         \\ 
         \&         \&         \&         \&         \& \nodef  \&         \&         \&         \&         \&         \\ 
         \&         \&         \&         \&         \&         \& \nodeg  \&         \&         \&         \&         \\
};

\path[ultra thick, black] (f) edge (g)
	(e) edge (f)
	(d) edge (e)
	(b) edge (c) edge (d)
	(a) edge (b);
\end{tikzpicture}}

}

%

\section{A multivariate version of the Catalan family}
\subsection{A generic factorisation of $\sigma_a$} 

The derivation of the Catalan idempotents presented in  \cite[Sec. 10]{MNT} can be interpreted as a Birkhoff factorisation
of the character of $QSym$ defined by 
\begin{equation}\label{eq:defphi}
\varphi(M_I)=\begin{cases}a^n&\text{if $I=(n)$}\\ 0&\text{otherwise,}\end{cases}
\end{equation}
for a certain choice of $a$ in a Rota-Baxter algebra of functions $\R$, the Rota-Baxter map being the the multiplpication by the indicatrix of $\R^+$.

Let us now look at the generic factorisation of this character, for an arbitrary Rota-Baxter algebra
$\AA=\AA_+\oplus \AA_-$, with projectors $P_+$ and $P_-$. Under the embedding \eqref{eq:SymInH} of $\Sym$
into $\Hh_{NCK}^*$, we have
\begin{equation}
\sigma_a = \sum_{n\ge 0}a^n S_n =\sum_{F\in\Ff}\varphi(Y_F)X_F.
\end{equation}
Writing the Birkhoff factorization $\varphi^+=\varphi^-\star\varphi$ as
\begin{equation}
	\sigma_a^+ =\sigma_a^-\sigma_a, \quad
\sigma_a^\pm =\sum_{F\in\Ff}\varphi^\pm(Y_F)X_F, 
\end{equation}
we can easily calculate $\varphi^\pm$ by remarking that
\begin{equation}
\lambda_{-a} = \sum_{n\ge 0}(-a)^nX_{\underbrace{\bullet\bullet\cdots\bullet}_{n}}=:\sum_{F\in\Ff}\alpha(Y_F)X_F
\end{equation}
where the character $\alpha$ is defined \textbf{on trees} by
\begin{equation}
\alpha(Y_T)=\begin{cases}-a&\text{if $T=\bullet$}\\0&\text{otherwise.}\end{cases}
\end{equation}
Since $\sigma_a^+\lambda_{-a}=\sigma_a^-$, we have $\varphi^+\star \alpha=\varphi^-$, which gives, for $T=B_+(F)$,
\begin{equation}
\varphi^+(Y_T)+\varphi^+(Y_F)\alpha(Y_\bullet) = \varphi^+(Y_T)-\varphi^+(Y_F)a=\varphi^-(Y_T)
\end{equation}
an applying $P_+$ and $P_-$, we obtain the recursive formulas
\begin{align}\label{recphi}
\varphi^+(Y_T) &= P_+(\varphi^+(Y_F)a),\\
\varphi^-(Y_T)&= -P_-(\varphi^+(Y_F)a).
\end{align}

\subsection{Example: polar part of  a Laurent series}
Let us now take $\AA= \C[z^{-1},z]]$ with $\AA_+=z^{-1}\C[z^{-1}]$
and $\AA_-=\C[[z]]$, and let
\begin{equation}
a = a(z) = \sum_{n\ge 0}a_nz^{n-1}.
\end{equation}
Here, $P_+(f)$ is defined as the polar part.

\medskip
{\footnotesize
We have, writing $a_{i_1\cdots i_r}$
for $a_{i_1}\cdots a_{i_r}$,
\begin{align*}
\varphi^+(Y_\bullet) &= P_+(a)=\frac{a_0}{z}  \\
\varphi^+(Y_{\bullet\bullet}) &=P_+(a)^2=\frac{a_{00}}{z^2}   \\
\varphi^+(Y_{\arbuga}) &= P_+\left(\frac{a_0}{z}\left(\frac{a_0}{z}+a_1+\cdots\right)\right) = \frac{a_{00}}{z^2}+\frac{a_{01}}{z}  \\
\varphi^+(Y_{\arbdgb}) &=  P_+\left(\frac{a_{00}}{z^2}\left(\frac{a_0}{z}+a_1+a_2z\cdots\right)\right)=      \frac{a_{000}}{z^3}+\frac{a_{010}}{z^2}+\frac{a_{002}}{z}    \\
\varphi^+(Y_{\arbdga}) &=  \frac{a_{000}}{z^3} + \frac{a_{010}}{z^2} + \frac{a_{001}}{z^2} + \frac{a_{011}}{z} + \frac{a_{002}}{z}.  
\end{align*}
}
On these examples, we can observe the following explicit description:
\begin{theorem}\label{th:phiX}
For any tree $T$, the value of the character $\varphi^+$ on a tree $Y_T$ is given by
 \begin{equation}
\varphi^+(Y_T)=\sum_{F\ge T}a_F z^{-r(F)}\label{eq:phi+tamari}
\end{equation}
where the sum is over the Tamari order on plane forests, $r(F)$ denotes the number of roots of $F$ and $a_F=a_{c_1}\cdots a_{c_n}$  
if $c_1\cdots c_n$ is the reverse Polish code of $F$.
\end{theorem}

\Proof This is an immediate consequence of the recursive description of Tamari intervals given below.
%
\qed

\subsection{A recursive description of the Tamari order}

\begin{theorem}
For a forest $F$, let $f(F)$ be the formal sum
\begin{equation}
f(F)=\sum_{G\ge F}G.
\end{equation}
Then for a tree $T=B_+(F)$,
\begin{equation}
f(T)= \sum_{F=F_1F_2}f(F_1)B_+(f(F_2)).
\end{equation}
\end{theorem}

In other words, the formal sum of the reverse Polish codes of the forests $G\ge B_+(F)$
is obtained by the following process: for each tree $T'=B_+(F')\ge T$, write down
the reverse Polish code $a_{F'}=a_{c_1}\cdots a_{c_{n-1}}$, and take the sum
of the  words $a_{F'}a_i$ for $i=0,\ldots,r$, where $r$ is the number of connected components
of $F'$. This amounts to encoding $F'$ by $\frac{a_{F'}}{z^r}$ and taking the polar part
of $\frac{a_{F'}}{z^r}a(z)$, which implies Theorem \ref{th:phiX}.

\bigskip
{\footnotesize
For example with $T= \arbtgc$, the codes of the trees $T'\ge T$ are
$0021, 0102, 0003$. The above process gives for each of them
\begin{align*}
0021 &\rightarrow 0020 + 0021,\\
0102 &\rightarrow 0100 + 0101 + 0102,\\
0003&\rightarrow 0000+0001+0002+0003,
\end{align*}
which are indeed the codes of the 9 forests  $G \ge T$.
}

\bigskip
\Proof
Recall the cover relation of the Tamari order on plane trees: starting from a
tree $T$ and a vertex $x$ that is neither its root or a leaf, the trees
$T'>T$ covering $T$ are obtained by cutting off the leftmost subtree of $x$
and grafting it back on the left of the parent of $x$.

So all elements in the Tamari order above a given element are obtained by a
sequence of such moves 
which can be encoded by a sequence of numbers
recording on which node each cut is done.

We shall prove the result for  a forest containing a single tree since the proof
works in the same way  with a general forest. 
We shall actually prove a stronger result: all
trees above a given tree can be obtained by a sequence of cuts where no cut is
done inside a subtree that was already cut.

To see that, number the internal nodes of a tree in prefix order, so that any
node has a label smaller than its descendants. Now, the path from a tree to a tree
above  it corresponds to a word on these labels  recording in which order the
nodes were cut. Assume that there is somewhere a factor $1i$ where $i>1$.
Then we shall see that this factor can be rewritten either as $i1$ if $i$ is
not the leftmost child of $1$ at this step or as $i11$ it it is.

First, if $i$ is the leftmost child of the root, applying $1$ then $i$
leads to a forest containing three trees: the left subtree of $i$, the
remaining part of the tree of root $i$ without its left child and the
remaining part of the whole tree without its left child.
One easily checks that we get the same result by applying $i11$ to the tree.

Moreover, if $i$ is not the left-most child of the root, then $1$ and
$i$ commute since they work in separate parts of the tree.

So by induction, any word sending a tree $T$ to a tree $T'$ below it can be
rewritten as a word where its 1s are at the end. The same applies to any
element of the tree, whence the result.
\qed
 
\subsection{Grouplikes and primitives in $\Sym$}

By definition, the series
\begin{equation}
\Cc:=\sigma_a^+|_{z=1} \quad \text{and}\quad \DD = \Res_{z=0}\sigma_a^+
\end{equation}
are in $\Sym$. We shall see that $\DD$ is a multiparameter version of the Catalan idempotent
of \cite{MNT,FMNT}, which is obtained by the choice $a(z)=\frac{a}{z}+\frac{b}{1-z}$. 

Indeed, recall that the basis $C_F$ of $\Hh_{NCK}^*$ is defined by 
\begin{equation}
C_F=\sum_{G\le F}X_G.
\end{equation}
Thus, 
\begin{equation}
\Cc = \sum_F\left(\sum_{G\ge F}a_G\right)X_F = \sum_G a_G\sum_{F\le G}X_F =\sum_G a_G C_G,
\end{equation}
where $a_G$ is the (commutative) product of the code of $G$, and since taking the residue amounts to
restricting the sum to trees,
\begin{equation}
\DD= \sum_T a_T C_T
\end{equation}
which gives back the expression obtained in \cite{FMNT} for $a(z)=\frac{a}{z}+\frac{b}{1-z}$. 

The possible values of $a_T$ correspond to partitions $\lambda$ of $n-1$. The sums
\begin{equation}
\DD_\lambda := \sum_{a_T=a_\lambda} C_T
\end{equation}
are therefore Lie quasi-idempotents of the descent algebra. 

\medskip
{\footnotesize
For example,
\begin{align*}
\DD_{(3)} &= C_{3000} = \bar\Psi_4,\\
\DD_{(21)} & = C_{2100}+C_{2010}+C_{1200} = R_4-R_{22}+R_{121}-R_{111}+\Psi_4+\bar\Psi_4,\\
\DD_{(111)} &= \Psi_4.
\end{align*}
}
\medskip

\subsection{Expansions in $\Sym$}

To compute the expansions of $\Cc$ and $\DD$ on the usual bases of $\Sym$, we start with
the Birkhoff recurrence \cite{Man}
\begin{align}
\varphi^-(x)&=-P_-\big(\varphi(x)+\sum_{(x)}\varphi^-(x')\varphi(x'')\big)\label{eq:recBirk1}\\
\varphi^+(x)&=P_+\big(\varphi(x)+\sum_{(x)}\varphi^-(x')\varphi(x'')\big)\label{eq:recBirk2}.
\end{align}
This gives immediately 
\begin{align}
\varphi^-(M_n)=-P_-(a^n),\ \varphi^+(M_n)=P_+(a^n),\\
\varphi^-(M_{ij})=P_-(P_-(a^i)a^j),\ \varphi^+(M_{ij})=-P_+(P_-(a^i)a^j),\ldots
\end{align}
and by induction,  we arrive at the proposition below.

Let $a$ an element of a Rota-Baxter algebra $\AA$.
We set $P_{\emptyset}^{ \emptyset}(a)=1_{\KK}$ and for  $I=(i_1,\dots,i_n)$, $\tmmathbf{\varepsilon} 
= (\varepsilon_1, \ldots, \varepsilon_n) \in \{+,-\}^n$,
\begin{equation}
P^{I}_{\tmmathbf{\varepsilon}}(a)=P_{\varepsilon_n}\left(P^{I'}_{\tmmathbf{\varepsilon}'}(a)a^{i_n}\right)
\end{equation}
where $I'=(i_1,\dots,i_{n-1})$ and $\tmmathbf{\varepsilon}'=(\varepsilon_1, \ldots, \varepsilon_{n-1})$.
We also write for short $P_{\tmmathbf{\varepsilon}}(a)=P_{\varepsilon_1, \ldots, \varepsilon_n}(a)$ 
the element of $\AA^{\varepsilon_n}$ equal to $P^{I}_{\tmmathbf{\varepsilon}}(a)$ where $i_1=i_2=\dots=i_n=1$.

\bigskip
{\footnotesize
For instance,\[
P^{1,2,3}_{+,-,-}=P_-(P_-(P_+(a)a^2)a^3)\text{ and } P_{+,+,-}(a)=P_-(P_+(P_+(a)a)a).
\]
}

\begin{proposition}\label{prop:basisexp} Let $a$ an element of a Rota-Baxter algebra $\AA$. Then,
\begin{align}
	\sigma_a^+&=\sum_I (-1)^{l(I)-1}P^{I}_{(-)^{l(I)-1},+}(a)S^I\\
	&=\sum_I (-1)^{|I|+l(I)}P^{I}_{(+)^{l(I)}}(a)\Lambda^I\\
	&=1+\sum_{\tmmathbf{\varepsilon}\in \mathcal{E}}P_{\tmmathbf{\varepsilon},+}(a)R_{\tmmathbf{\varepsilon},\bullet}
\end{align}
and
\begin{align}
	\sigma_a^-&=\sum_I (-1)^{l(I)}P^{I}_{(-)^{l(I)}}(a)S^I\\
	&=\sum_I (-1)^{|I|+l(I)-1}P^{I}_{(+)^{l(I)-1},-}(a)\Lambda^I\\
	&=1-\sum_{\tmmathbf{\varepsilon}\in \mathcal{E}}P_{\tmmathbf{\varepsilon},-}(a)R_{\tmmathbf{\varepsilon},\bullet}
\end{align}
\end{proposition}

We use in the last equation the {\it signed ribbon basis} of $\Sym$ (see \cite{MNT}),
 which is a slight modification
of the noncommutative ribbon Schur functions: for any sequence of signs $\mathbf{\varepsilon}=(\varepsilon_1,\dots,\varepsilon_{n-1}$
\begin{equation}
R_{\mathbf{\varepsilon} \bullet}=(-1)^{l(I)-1}R_I\quad
(R_\emptyset=1,\ R_{\bullet}=R_1)
\end{equation}
where $I=(i_1,\dots,i_r)$ is the composition of $n$ such that
\begin{equation}
D(I):=\{i_1,i_1+i_2,\dots,i_1+\dots+i_{r-1}\}=\{1\leq i \leq n-1 \ ;\ \varepsilon_i=-\}.
\end{equation}
\Proof The expansions on the $S^I$ follow immediately from the recurrence \eqref{eq:recBirk1}-\eqref{eq:recBirk2}. The other ones admit an interesting explanation
in terms of the free Rota-Baxter algebra on one generator, which can be realized a subalgebra of the algebra of sequences
of multivariate polynomials, with pointwise addition and product.

Let $\xx$ be a sequence of variables $\xx = (x_1,x_2,x_3,\ldots)$
and for a sequence $\zz$, define
\begin{equation}
 R(\zz) = (0,z_1,z_1+z_2,z_1+z_2+z_3,\ldots)
\end{equation}
This is a Rota-Baxter operator of weight $1$. It generates from $\xx$ the free Rota-Baxter algebra
 $\mathfrak{A}(\xx)$ \cite{Rota}.
Define
\begin{equation}
 P_+ =-R,\quad P_- = I-P_+,
\end{equation}
which are now of weight $-1$.
Set also $V(\zz)=(z_2,z_3,\ldots)$.

The subalgebra generated by the $R(\xx^n)$ is isomorphic to $Sym$:
\begin{equation}
R(\xx^n) = (p_n(0),p_n(x_1),p_n(x_1,x_2),\ldots):=\tilde p_n
\end{equation}
but there is also an embedding  of $QSym$ in  $\mathfrak{A}(\xx)_+:=P_+(\mathfrak{A}(\xx))$  given by the same rule 
\begin{equation}
f\mapsto \tilde f :=(f(0),f(x_1), f(x_1,x_2),\ldots)
\end{equation}
Its image under $V$ gives an embedding of $QSym_+$ into $\mathfrak{A}(\xx)_-$.
 
It is now easy to show by induction that 
for $I=(i_1,\ldots,i_r)$, 
\begin{equation}
	\tilde M_I = R(\tilde M_{i_1\cdots i_{r-1}}\xx^{i_r}).
\end{equation}
If we number the applications of $R$ in the above expression of $M_I$, 

\medskip
{\footnotesize
\noindent for example
\begin{equation}
\tilde M_{i_1i_2i_3i_4}=R_4(R_3(R_2(R_1(a^{i_1})\xx^{i_2})\xx^{i_3})\xx^{i_4})
\end{equation}
}

\medskip
\noindent and replace some $R_j$ by $R'=I+R=P_-$, the result is now $\tilde M_I+\tilde M_{I'}$,
where $I'=(i_1,\ldots,i_j+i_{j+1},i_{j+2},\ldots)$. 

\medskip
{\footnotesize
For example, 
\begin{equation}
\tilde M_{i_1i_2i_3i_4}\mapsto R_4(R_3((I+R_2)(R_1(\xx^{i_1})\xx^{i_2})\xx^{i_3})\xx^{i_4})
	=\tilde M_{i_1i_2i_3i_4}+\tilde M_{i_1,i_2+i_3,i_4}.
\end{equation}
}

\medskip
Applying this to $M_{1^n}$, the replacement of $R_i$ by $R'_i$ has the effect of removing $i$ from the descent set
of $1^n$. Iterating, we see that replacing $R_{d_1},\ldots,R_{d_k}$ by $R'$ yields all compositions
whose descent set contains the complement of $\{d_1,\ldots,d_k\}$. The result is therefore $\tilde F_{\bar I^\sim}$

Let us now look at the coefficient of $S^I$ in $\sigma_\xx^+$. For example, for $I=(i,j,k)$, this is
\begin{align*}
	(-1)^3P^{ijk}_{--+}(\xx)&= R(((I+R)((I+R)\xx^i)\xx^j)\xx^k)\\
	&=R(\xx^{i+j+k}+R(R(\xx^i)\xx^{j+k}+R(\xx^{i+j})\xx^k+R(R(\xx^i)\xx^j)\xx^k)\\
	&= \tilde M_{i+j+k}+\tilde M_{i,j+k}+\tilde M_{i+j,k}+\tilde M_{ijk}\\
	&= \sum_{J\le I}\tilde M_J\\
	&=(-1)^3 \widetilde{M_{ijk}(-X)}.
\end{align*}
Thus, the Birkhoff factorisation of the character \eqref{eq:defphi} of $QSym$ is given by

\begin{equation}
	\varphi^+(M_I)=\widetilde{M_I(-X)},\quad  
	\varphi^-(M_I)=V\widetilde{M_I(-X)},
\end{equation}
and
\begin{equation}
 \varphi^+(F_I)=\tilde F_I(-X)=(-1)^{|I|}\tilde F_{\bar I^\sim},\quad \varphi^-(F_I) =V\varphi_+(\tilde F_I).
\end{equation}
Now, the coefficient of $\Lambda^I$ in $\sigma_\xx^+$ is $\tilde M_I(X)$, whence the second equalities in Prop. \ref{prop:basisexp},
and that of $R_I$ is, up to sign, $\tilde F_{\bar I^\sim}(X)$, which can be expressed as 
$\pm P_{\tmmathbf{\varepsilon},+}(\xx)$.
\qed

\subsection{Combinatorial interpretation of the coefficients}

Evaluating the above expression for $\varphi^+(M_I)$ on $\xx=a(z)$, but without assuming that the $a_i$ commute
yields a set of words which can be characterized by certain inequalities involving partial sums
of the subscripts. Recall that
\begin{equation}
\sigma_a^+ = \sum_{I} \varphi^+(M_I) S^I,
\end{equation}
where $\varphi^+(M_n)=P_+(a^n)$, and for $I=(I',i_p)$.
\begin{equation}
\left\{
\begin{array}{l}
\varphi^+(M_I) =   P_+ (\varphi^-(M_{I'}) a^{i_p})   \\
\varphi^-(M_I) = - P_- (\varphi^-(M_{I'}) a^{i_p}).
\end{array}
\right.
\end{equation}
So if  $I=(i_1,\dots,i_p)$, the evaluation of $\varphi^+(M_I)$
amounts to computing $P_-(M_{i_1})$, then multiply by
$a^{i_2}$, then apply $P_-$, then multiply the result by $a^{i_3}$, and so on,
up to the last step where one applies a $P_+$ instead of $P_-$.

Thus, up to a global sign $(-1)^{\ell(I)-1}$, $\varphi^+(M_I)$ is a sum of monomials
in $z^{-1}$ and the $a_i$.
Such a monomial is a product
of $n$ terms of the series $a$ which survive the sequence of  $P_-$ and the final
$P_+$. Writing this product as a word, considering that the $a_i$ do not
commute, and replacing $a_i$ with $i$ (ignoring the power of $z$ that can be
reconstituted in the end),  we  obtain a word $w=w_1\dots w_n$ 
over the integers such that
\begin{equation}
\left\{
\begin{array}{l}
w_1+\dots + w_{d_k} \geq d_k, \text{\ for all $k<p$}, \\
w_1+\dots + w_{i_1+\dots+i_p} < i_1+\dots+i_p,
\end{array}
\right.
\end{equation}
where $\{d_1=i_1,d_2=i_1+i_2,\dots,d_{p-1}=i_1+\dots+i_{p-1}\}$ is the
\emph{descent set} $D(I)$ of $I$. 
Denote this set of words by $S(I)$.

\medskip
{\footnotesize
Let us check this observation on $\varphi^+(M_{112})$, which is
\begin{equation}
\frac{a_3 a_0^3}{z} + \frac{4a_2a_1a_0^2}{z} + \frac{2a_1^3a_0}{z}
+ \frac{a_2a_0^3}{z^2} + \frac{a_1^2a_0^2}{z^2}.
\end{equation}
It is indeed obtained from the $9$ words $w=w_1w_2w_3w_4$ satisfying
\begin{equation}
w_1\geq 1, \quad
w_1+w_2\geq 2, \quad
w_1+\dots+w_4 < 4,
\end{equation}
that are
\begin{equation}
3000,\ 2100,\ 2010,\ 2001,\ 1200,\ 1110,\ 1101,\ 2000,\ 1100,
\end{equation}
by sending each value $i$ to $a_i z^{i-1}$.
}

\medskip
For a word over the integers, define 
\begin{equation}
w_{1:k}:=\sum_{i=1}^kw_i
\end{equation}
and let
\begin{equation}
W(I)=\{w|w_{1:k}\ge k \ \text{if $k\in D(I)$ and}\ w_{1:k}<k\ \text{otherwise}\},
\end{equation}
so that
\begin{equation}
S(I)=\bigsqcup_{J\ge I} W(J).
\end{equation}
Thus,
if one writes
as an intermediate expression for $\varphi^+(M_I)$ the sum
\begin{equation}
\begin{split}
\sum_I  (-1)^{\ell(I)-1} S^I \sum_{w\in S(I)} w
&= 
\sum_I  (-1)^{\ell(I)-1 }S^I \sum_{J\ge I}\sum_{w\in W(J)} w\\
&=\sum_J\sum_{w\in W(J)} w\sum_{I\le J}(-1)^{\ell(I)-1}S^I
\end{split}
\end{equation}
one can see that the coefficient of a word $w\in W(J)$ is, up to a sign $(-1)^{\ell(J)-1}$, the ribbon $R_J$.

So the expansion of $\sigma_a^+$ in the ribbon basis  is  obtained  by listing
the words $w=w_1\dots w_n$ satisfying $w_1+\dots+w_n<n$ (counted by
the binomial $\binom{2n-1}{n}$). Each such $w$ belongs to a unique $W(I)$,
which determines its  coefficient $(-1)^{\ell(I)-1}R_I$, and a
factor  $z^{w_{1:n}-n}$  

\medskip
{\footnotesize
For example, here are all possible words for $n=3$ with the corresponding
compositions:
\begin{equation}
\begin{array}{llllllllll}
000 & 001 & 010 & 100 & 002 & 020 & 200 & 011 & 101 & 110 \\
 3  &  3  &  3  &  12 &  3  &  21 & 111 &  3  &  12 & 111
\end{array}
\end{equation}

For $n=4$, here is the complete list of all words contributing to each $R_I$:
\begin{equation}
\begin{array}{ll}
4    & 0000,\ 0100,\ 0010,\ 0001,\ 0110,\ 0101,\ 0020,\ 0011,\ 0002,\\
     & 0111,\ 0102,\ 0021,\ 0012,\ 0003, \\
31   & 0120,\ 0030, \\
22   & 0200,\ 0201, \\
13   & 1000,\ 1010,\ 1001,\ 1011,\ 1002, \\
211  & 0300,\ 0210, \\
121  & 1020, \\
112  & 2000,\ 1100,\ 2001,\ 1101, \\
1111 & 3000,\ 2100,\ 2010,\ 1200,\ 1110.
\end{array}
\end{equation}
}

\medskip
We already know from previous works \cite{MNT,FMNT} that if $a_0=a, a_i=b$ for $i>0$, the coefficient
of a $R_I$ is (up to a global sign) a product of Narayana polynomials. Since the coefficients
in the general case are sums of monomials with the same sign, this implies that the cardinalities
of the sets $W(I)$ are products of Catalan numbers. This can be seen directly as follows.

Recall the correspondence between Łukasiewicz words (Polish codes of plane trees) and
Dyck paths.
The code of a plane tree is obtained by labelling each node by the number of its descendants, and traversing it in prefix order.

\medskip
{\footnotesize
An example would be
\begin{equation}
w= 40201200010
\end{equation}
}

\medskip
These codes are characterized by the following property: if one forms a word $u$ by subtracting 1 to each entry of $w$, the partial sums $u_{1:i}$ are all nonnegative, except for
the last one which is $-1$. 

\medskip
{\footnotesize
On our example, 
\setcounter{MaxMatrixCols}{20}
\begin{equation}
\begin{matrix}
4 & 0 & 2 & 0 & 1 & 2 & 0 & 0 & 0 & 1 & 0\\
3 &-1 & 1 &-1 & 0 & 1 &-1 &-1 &-1 & 0 &-1\\
3 & 2 & 3 & 2 & 2 & 3 & 2 & 1 & 0 & 0 & -1
\end{matrix}
\end{equation}
}

\medskip
This characterization means that if one replaces each integer $i$ by the word
$a^ib$, one obtains a word $wb$, where $w$ is a Dyck word\footnote{Here, the letter $a$ stands for an upstep and $b$ for a downstep.}. 

\medskip
{\footnotesize
On our exemple, this yields
\begin{equation}
aaaab.b.aab.b.ab.aab.b.b.b.ab\cdot b
\end{equation}
}

\medskip
The partial sums $u_{1:i}$ give the height of the corresponding Dyck path after the $i$th $b$.

This description can be extended to the sets $W(I)$. 
The word obtained by replacing each entry $k$ by $a^kb$ in $w$ encodes a lattice path starting
at the origin, and ending at $(2n+1,-1)$. 
Applying the transformation $u_i=w_i-1$
to $W(I)$ results into the set of words
\begin{equation}
U(I)=
\{u|u_{1:k}\ge 0 \ \text{if $k\in D(I)$ and}\ u_{1:k}<0\ \text{otherwise}\}.
\end{equation}
Again, the partial sums $u_{1:i}$ of such words record the heights attained
by the lattice path associated with $w$ after the $i$th $b$.

Represent a composition $I=(i_1,\dots,i_p)$  of $n$ as a
sequence of $n$ symbols $+$ and $-$ with a $-$ in position $k$ if
$k$ is a descent of $I$, and a $+$ otherwise.

\medskip
{\footnotesize
For example,  $312$ is represented as  $++--++$ and $3111$ as
$++---+$. 
}

\medskip
Then, the cardinality of $W(I)$
is $\prod_i C_{i}$ where $i$ runs over the lenghts of blocks of identical signs.

\medskip
{\footnotesize
For example,  $W(312)$ contains $C_2^3=8$ words and $W(3111)$ has $C_2C_3=10$ elements.
}

\medskip
Indeed, the blocks of symbols $+$ correspond to sections of the path associated with $w$ lying
under the horizontal axis, and the blocks of $-$ to sections where it remains above the axis.
The sections of the path determined by these blocks are alternatively Dyck paths or negative
of Dyck paths, whence the product of Catalan numbers. Counting them by number of peaks gives back
the products of Narayana polynomials already mentioned.

\medskip
{\footnotesize
For example, let us decompose  $W(4111)$.
The corresponding signed word is $+++---+$. 
There should be $25$ such words. Let us write these as a $5\times5$ square
where words on the same column have same first three values.
\begin{equation}
\begin{array}{lllll}
0006000 & 0015000 & 0105000 & 0024000 & 0114000 \\
0005100 & 0014100 & 0104100 & 0023100 & 0113100 \\
0005010 & 0014010 & 0104010 & 0023010 & 0113010 \\
0004200 & 0013200 & 0103200 & 0022200 & 0112200 \\
0004110 & 0013110 & 0103100 & 0022110 & 0112110
\end{array}
\end{equation}
The path corresponding to $0004200$ is $bbbaaa.abaabb.b$,
and that corresponding to $0112200$ is $bababa.abaabb.b$.
One can check that all pairs of Dyck paths are obtained. 
Note that in each row, the values $(w_4,w_5,w_6)$ are the
same if one replaces the fourth one by $w_4+(w_1+w_2+w_3)-3$.
The sequence of these values becomes
\begin{equation}
300, 210, 201, 120, 111,
\end{equation}
which is indeed the set of the first three values associated with the
composition $1111$, and the Polish codes of plane trees with 4 nodes except for their final 0.
}

\section{Lie idempotents of the descent algebra}

We shall now decribe the expansions of several Lie idempotents of the descent algebra on the $X$-basis.
To this aim, we shall need several versions of the $(1-q)$-transform.

Recall that on ordinary symmetric function, the alphabet $\frac{X}{1-q}$ is the set $\{q^ix_j|i\ge 0, x_j\in X\}$.
It can be extended to noncommutative symmetric functions by choosing a total order of the products $q^ia_j$, which can of
course be done in an infinity of ways, but only four of them are natural: take the lexicographic order on the pairs
$(q^i,a_j)$ or $(a_j,q^i)$, keeping the original order on $A$ and ordering the $q^i$ in ascending or descending order
of the exponents. This leads to four possible definitions of the $(1-q)$-transform as the respective inverses of the
above transforms. In the sequel we shall define them directly by specifying the image of the $S_n$.

\subsection{Dynkin}
\begin{proposition}
The right Dynkin $\bar\Psi_n=[1,[2,[3,\ldots[n-1,n]\ldots]]]$
is the sum of all trees
\begin{equation}
\bar\Psi_n = \sum_{|T|=n}X_T.
\end{equation}
and the left Dynkin $\Psi_n=[\ldots[[1,2],3],\ldots,n]$ is the linear tree $X_{L_n}$
\begin{equation}
\bar\Psi_n = ((X_\bullet\triangleright X_\bullet)\cdots)\triangleright X_\bullet.
\end{equation}
\end{proposition}

\Proof
We first apply Theorem \ref{th:coeffX} to
$X=1-q$, defined by
\begin{equation}
S_n((1-q)A)= (1-q)\sum_{k=0}^n(-q)^k R_{1^k,n-k}(A),
\end{equation}
so that $\Psi_n(A) = \frac1{1-q}S_n((1-q)A)|_{q=1}$, and
$F_I(1-q)$ is nonzero only for  $I$ of the type $(1^k,n-k)$. 

Every forest with $k+1$ leaves has a unique maximal linear extension of this shape,
obtained by reading its leaves from right to left and then taking the postorder reading
of the remaining nodes. It has therefore ${k\choose i}$ linear extensions of
shape $(1^i,n-i)$ for $0\le i\le k$, so that $\Gamma_F(1-q)=(1-q)^{k}$ is divisible
by $(1-q)^2$ except for $k=1$, which means that $F=L_n$ is a linear tree.

To deal with $\bar \Psi_n$, we need another version of the $1-q$ transform, denoted by
$1+(-q)$, and défined\footnote{
	This strange notation is justified by the fact  that addition of alphabets is not commutative, and that $X-Y$ is defined as $(-Y)+X$, {\it cf.} \cite{NCSF2}.} 
by $F_I(1+(-q))=(1-q)(-q)^k$ if
$I=(n-k,1^k)$ and 0 otherwise,  
so that $\bar \Psi_n(A) = \frac1{1-q}S_n((1+(-q))A)|_{q=1}$.
A permutation of shape $I=(n-k,1^k)$ cannot be a linear extension of a tree, unless $k=0$, in which case it is the identity,
the common linear extension of all trees. Thus, $\bar\Psi_n$ is the sum of all trees with $n$ nodes.
\qed

\subsection{Eulerian idempotents}
Take the binomial alphabet $\alpha$ defined by $\sigma_1(\alpha A)=\sigma_1^\alpha$,
so that $M_I(\alpha)={\alpha\choose\ell(I)}$,
and $F_I(\alpha)={\alpha+n-r\choose n}$ where $n=|I|$ and $r=\ell(I)$.
Then, the Solomon idempotent 
$\varphi$ (often denoted by  $\Omega$, and also known as the first Eulerian idempotent) is given by
\begin{equation} 
\varphi := \log \sigma_1 = \left.\frac{d}{d\alpha}\right|_{\alpha=0}\exp{\alpha\varphi} =\left.\frac{d}{d\alpha}\right|_{\alpha=0}\sigma_1(\alpha A),
\end{equation}
so that the coefficient of  $X_T$ in $\varphi$ is
\begin{equation} 
\left.\frac{d}{d\alpha}\right|_{\alpha=0}\Gamma_T(\alpha).
\end{equation}
Equivalently, with the notation of Theorem \ref{th:coeffX}
\begin{equation}
\sum_F \chi_F(\alpha)X_F = \lambda_1(A)^\alpha
=\exp\left\{\alpha\sum_{n\ge 1}(-1)^{n-1}\varphi_n\right\}
\end{equation}
and for a forest of degree $n$,
\begin{equation}
\left.\frac{d}{d\alpha}\right|_{\alpha=0}\chi_F(\alpha) = (-1)^{n-1}( Y_F,\varphi_n)
\end{equation}
so that
\begin{equation}
\varphi_n = (-1)^{n-1}\left.\frac{d}{d\alpha}\right|_{\alpha=0}
\sum_{|T|=n} \chi_T(\alpha)X_T
\end{equation}
which contains only trees, since
$\varphi$ is a Lie series.

The polynomial
$\chi_T(t)$ is the evaluation of the tree $T$ obtained by putting $t$ in each leaf,
the operator ``discrete integral of the product of the subtrees'' 
\begin{equation}
\Sigma:\ t^p\mapsto \Sigma_0^t s^p\delta s = \frac{B_{p+1}(t)-B_{p+1}(0)}{p+1}
\end{equation}
in each internal node, and multplying the result by
$(-1)^{n-1}$ (the $B_k$ 
are the Bernoulli polynomials).

Indeed, if $T=B_+(T_1\cdots T_k)$,  $\chi_T(t)$ satisfies the difference equation
\begin{equation}
	\Delta \chi_T(t) =\chi_{T_1}(t)\cdots \chi_{T_k}(t)
\end{equation}
which can bee seen as follows. First, $\chi_T(t)=\<Y_T,\lambda_1^t\>$, so that
\begin{align*}
	\Delta\chi_T(t) &= \<Y_T,\lambda_1^t(\lambda_1-1)\> =\<\Delta Y_T,\lambda_1^t\otimes(\lambda_1-1)\>\\
	&= \sum_{(T)}\<Y_{T(1)}\otimes Y_{T(2)},\lambda_1^t\otimes(X_\bullet+X_{\bullet\bullet}+\cdots)\>\\
	&= \<Y_{T_1}\cdots Y_{T_k},\lambda_1^t\> \quad\text{since the only nonzero term is obtained for $T(2)=\bullet$}\\
	&= \chi_{T_1}(t)\cdots \chi_{T_k}(t).
\end{align*}
This formula has been first obtained in \cite{WZ} by a more complicated argument.

The  coefficients of the polynomial
$(-1)^{|T|}\chi_T(t)$ given the expansion of the other Eulerian idempotents is the forests basis.
This is equivalent to the description of the ``formal flow'' given in
\cite{WZ}.
The coefficient of $\alpha^k$ in $\sigma_1^\alpha$ is
\begin{equation}
\frac1{k!}\sum_{\ell(I)=k} \varphi^I
\end{equation}
hence the coefficient of $X_F$ in
\begin{equation}
e_n^{(k)}= \frac1{k!}\sum_{I\vDash n}\varphi^I
\end{equation}
is ({\it cf.} Eq. \eqref{eq:chi})
\begin{equation}
[\alpha^k]\Gamma_F(\alpha)=
(-1)^{|F|}[\alpha^k]\chi_F(\alpha).
\end{equation}

\subsection{$q$-idempotents and a two-parameter series}

In \cite{NCSF2}, it is proved that, for the usual definition of $\frac{A}{1-q}$
\begin{equation}
	\varphi_n(q) = \frac{1-q^n}{n}\Psi_n\left(\frac{A}{1-q}\right) = 
	{1 \over n} \ \sum_{|I|=n} \
{ (-1)^{\ell (I)-1} \over \qbin{n-1}{\ell (I) -1}} \
q^{\maj(I) - \ssbin{\ell (I)}{2}} \ R_I(A)
\end{equation}
is a Lie idempotent, interpolating between the Solomon idempotent $\varphi_n$ (for $q=1)$, 
the two Dynkin (for $q=0,\infty$) and Klyachko ($q=e^{2i\pi/n}$). Its expansion on the
preLie basis $x_\tau$ (hence also on $X_T$) is obtained by Chapoton in \cite{Cha1}.

One way to recover this result is to apply Theorem \ref{th:coeffX} to
the virtual alphabet
\begin{equation}
\frac{1-qt|}{|1-q}= (1-qt)\times\frac1{1-q}
\end{equation}
defined by \cite{NCSF2}
\begin{equation}
	S_n\left(\frac{1-qt|}{|1-q}A\right)=(1-qt)\sum_{k=0}^n(-qt)^k R_{1^k,n-k}\left(\frac{A}{1-q}\right) 
\end{equation}
so that
\begin{equation}
\Psi_n\left(\frac{A}{1-q}\right)= \frac1{1-qt}\left.S_n\left(\frac{1-qt|}{|1-q}A\right)\right|_{t=\frac1q}.
\end{equation}

The series denoted by $\sqx$ in \cite{Cha3} is essentially $\sigma_1\left(\frac{1-qt|}{|1-q}A\right)$.
Actually, Chapoton takes the opposite order on the alphabet of powers of $q$, and to recover the same
coefficients, we have to define $\sqx$ as
\begin{equation}\label{eq:defsqx}
	\sqx = \sigma_1(X_{q,t}A) := \prod_{i\ge 0}^\rightarrow\sigma_{q^i}(A)\prod_{j\ge 0}^\leftarrow\lambda_{-q^jt}(A).
\end{equation}
The functional equation satisfied by $f(t) := \sigma_1(X_{q,t}A)$ is then
\begin{equation}
	f(qt)=f(t)\sigma_{qt}(A)
\end{equation}
which is equivalent to \cite[(8)]{Cha3} after setting $t=1+(q-1)x$.

The coefficient of $\frac\tau{|{\rm Aut}(\tau)|}$ in 
$\sqx$ is thus obtained by setting $t=1+(q-1)x$ in $\Gamma_T\left(X_{q,t}\right)$.

For example, with  $T=10$, $\Gamma_T(A)=\M_{12}+\M_{11}$, hence $\Gamma_T(X)=M_2+M_{11}=h_2$ is a symmetric function,
and 
\begin{equation}
h_2\left(\frac{1-qt}{1-q}\right)=\frac{(1-qt)(1-q^2t)}{(1-q)(1-q^2)}=\frac{(1+qx)(1+q+q^2x)}{1+q}.
\end{equation}
Dividing by $1+qx$, and setting  $x=-1/q$, one  finds $\frac1{1+q}$, which is indeed the coefficient
of $X_{10}$ in the series $\bar\Omega_q$ defined in \cite[(45)]{Cha3}.

\subsection{Examples}

One can easily compute
$\Gamma_T(A)$ by the recurrence (obvious from the definition in terms of linear extensions)
\begin{equation}
\Gamma_{B_+(T_1\cdots T_k)}(A) = B(\Gamma_{T_1}\cdots\Gamma_{T_k}),
\end{equation}
where $B(\F_{\sigma}):=\F_{\sigma n}=\F_\sigma\succ\F_1$  ($n=|T_1|+\cdots+|T_k|+1$),
which yields by projection onto $QSym$
\begin{equation}
\Gamma_{B_+(T_1\cdots T_k)}(X) = B(\Gamma_{T_1}\cdots\Gamma_{T_k}),\quad\text{where $B(F_{i_1i_2\ldots i_r}):=F_{i_1,i_2\ldots,i_r+1}$}
\end{equation}

{\footnotesize
For example,
\begin{align}
\Gamma_{\arbuga}(X) &= F_2 \rightarrow {\alpha+1\choose 2}\\
\Gamma_{\arbdga}(X) &= F_3 \rightarrow {\alpha+2\choose 3}\\
\Gamma_{\arbdgb}(X) &= F_{12}+F_3 \rightarrow {\alpha+2\choose 3}+{\alpha+1\choose 3}\\
\Gamma_{\arbtge}(X) &= F_4 \rightarrow {\alpha+3\choose 4}               \\
\Gamma_{\arbtgc}(X) &=  F_{13} + F_4  \rightarrow {\alpha+3\choose 4}+  {\alpha+2\choose 4}\\
\Gamma_{\arbtgd}(X) &=  F_{22} + F_{13} + F_4     \rightarrow {\alpha+3\choose 4}+  2{\alpha+2\choose 4}      \\
\Gamma_{\arbtgb}(X) &=  F_{22} + F_{13} + F_4 \rightarrow {\alpha+3\choose 4}+  2{\alpha+2\choose 4}\\
\Gamma_{\arbtga}(X) &=  F_{112} + 2F_{22} + 2F_{13} + F_4\rightarrow {\alpha+3\choose 4}+  4{\alpha+2\choose 4}+{\alpha+1\choose 4}
\end{align}
which gives for the Eulerian idempotents
\begin{align}
e_4^{(1)}&= \frac1{4!}\left(6 X_{1110} + 4 X_{1200} + 2 X_{2010} + 2 X_{2100}  \right)=\varphi_4\\
e_4^{(2)}&= \frac1{4!}\left(9 X_{2100} + 6 X_{1010} + 6 X_{3000} + 10 X_{1200} + 9 X_{2010}\right. \nonumber\\
         &+ \left. 4 X_{2000} + 4 X_{0200} + 8 X_{1100} + 11 X_{1110} + 8 X_{0110}  \right)\\
e_4^{(3)}&= \frac1{4!}\left( 10 X_{2010} + 12 X_{0200} + 6 X_{1110} + 12 X_{0010} + 8 X_{1200} + 12 X_{0110}\right. \nonumber\\
         &+ \left. 12 X_{3000} + 12 X_{1100} + 12 X_{2000} + 12 X_{1010} + 12 X_{0100} + 10 X_{2100} + 12 X_{1000}  \right)\\
e_4^{(4)}&= \frac1{4!} \left( 3 X_{2100} + 8 X_{0200} + 8 X_{2000} + 2 X_{1200} + 4 X_{0110} + 4 X_{1100}\right.\nonumber\\
         &\left. + 3 X_{2010} + 12 X_{1000} + 12 X_{0100} + 6 X_{1010} + 6 X_{3000} + 24 X_{0000} + 12 X_{0010} + X_{1110} \right)
\end{align}
}

To recover Chapoton's coefficients for the two-parameter series, one has to use the other version of the $X$-basis,
defined by duality with the opposite coproduct on $\Hh_{NCK}$. This amounts to replacing $\Gamma(X)$ by $\Gamma'(X)=\omega(\Gamma(X))$,
that is, 
\begin{equation}
\Gamma_T(X_{q,t})=\omega(\Gamma_T)\left(\frac{1-qt|}{|1-q} \right).
\end{equation}
{\footnotesize
\begin{align}
\Gamma'_{\arbuga}\left(\frac{1-qt|}{|1-q} \right) &=\frac{{\left(q^{2} x + q + 1\right)} {\left(q x + 1\right)}}{q + 1} \\
\Gamma'_{\arbdga}\left(\frac{1-qt|}{|1-q} \right) &=\frac{{\left(q^{3} x + q^{2} + q + 1\right)} {\left(q^{2} x + q + 1\right)} {\left(q x + 1\right)}}{{\left(q^{2} + q + 1\right)} {\left(q + 1\right)}} \\
\Gamma'_{\arbdgb}\left(\frac{1-qt|}{|1-q} \right) &= \frac{{\left(q^{3} x + q^{2} x + q^{2} + q + 1\right)} {\left(q^{2} x + q + 1\right)} {\left(q x + 1\right)}}{{\left(q^{2} + q + 1\right)} {\left(q + 1\right)}}\\
\Gamma'_{\arbtge}\left(\frac{1-qt|}{|1-q}
\right) &=  \frac{{\left(q^{4} x + q^{3} + q^{2} + q + 1\right)} {\left(q^{3} x + q^{2} + q + 1\right)} {\left(q^{2} x + q + 1\right)} {\left(q x + 1\right)}}{{\left(q^{2} + q + 1\right)} {\left(q^{2} + 1\right)} {\left(q + 1\right)}^{2}}             \\
\Gamma'_{\arbtgc}\left(\frac{1-qt|}{|1-q} \right) &=  \frac{{\left(q^{3} x + q^{2} + q + 1\right)} {\left(q^{3} x + q^{2} + 1\right)} {\left(q^{2} x + q + 1\right)} {\left(q x + 1\right)}}{{\left(q^{2} + q + 1\right)} {\left(q^{2} + 1\right)} {\left(q + 1\right)}}\\
\Gamma'_{\arbtgd}\left(\frac{1-qt|}{|1-q} \right) &=      
\frac{{\left(q^{4} x + q^{3} x + q^{3} + q^{2} x + q^{2} + q + 1\right)} {\left(q^{3} x + q^{2} + q + 1\right)} {\left(q^{2} x + q + 1\right)} {\left(q x + 1\right)}}{{\left(q^{2} + q + 1\right)} {\left(q^{2} + 1\right)} {\left(q + 1\right)}^{2}}    \\
\end{align}

\begin{align}
\Gamma'_{\arbtgb}\left(\frac{1-qt|}{|1-q} \right) &=   
\frac{{\left(q^{4} x + q^{3} x + q^{3} + q^{2} x + q^{2} + q + 1\right)} {\left(q^{3} x + q^{2} + q + 1\right)} {\left(q^{2} x + q + 1\right)} {\left(q x + 1\right)}}{{\left(q^{2} + q + 1\right)} {\left(q^{2} + 1\right)} {\left(q + 1\right)}^{2}} \\
\Gamma'_{\arbtga}\left(\frac{1-qt|}{|1-q} \right) &= \frac{{\left(q^{6} x^{2} + q^{5} x^{2} + 2 \, q^{5} x + q^{4} x^{2} + 2 \, q^{4} x + q^{4} + 3 \, q^{3} x + q^{3} + 2 \, q^{2} x + 2 \, q^{2} + q + 1\right)} {\left(q^{2} x + q + 1\right)} {\left(q x + 1\right)}}{{\left(q^{2} + q + 1\right)} {\left(q^{2} + 1\right)} {\left(q + 1\right)}} 
\end{align}
}

\subsection{Appendix: noncommutative Ehrhart polynomials} 
In the introduction of \cite{Cha3}, Chapoton mentions that the coefficients of the series $\sqx$ are
$q$-analogues of Ehrhart polynomials (according to his definition given in \cite{Cha4}). These are actually
specialisations of the noncommutative Ehrhart polynomals, 
which  are defined only for the order polytopes of posets on $[n]$ \cite{BK,Wh}.

Recall the definition of the free generating function of a poset $P$
\begin{equation}
\Gamma_P(A) = \sum_{\sigma\in L(P)}\F_\sigma \in\FQSym
\end{equation}
where 
$L(P)\subseteq\SG_n$ is the set of linear extensions of $P$.
It is a morphism from the Malvenuto-Reutenauer Hopf algebra of special posets
towards $\FQSym$. In the sequel, we will only consider posets satisfying
$i<_P j\Rightarrow i<j$.

The order polytope $Q_P$ of $P$ is defined by the inequalities
$0\le x_i\le 1$ for $i\in P$ and $i<_Pj\Rightarrow x_i\le x_j$.

The Ehrhart polynomial $E_Q(t)$ computes the number of integral points of $nQ$ for $t=n$.
Moreover, $(-1)^nE_T(-n)$ is the number of interior integral points.

Since $nQ_P$ is the intersection of the cone $C_P$ defined by
$x_i\ge 0$ and $i<_Pj\Rightarrow x_i\le x_j$, and of a hypercube, one can form in $\WQSym$
the sum of the packed words of its integer points. The noncommutative Ehrhart polynomial of $Q_P$ is 
\begin{equation}
\sum_{u\in C(P)}\M_u = \Gamma_P(A)
\end{equation}
where $C(P)$ is the set of packed words satisfying $i<_Pj\Rightarrow u_i\le u_j$,
if one embeds
$\FQSym$ into $\WQSym$ by
\begin{equation}
\G_\sigma(A) = \sum_{\std(u)=\sigma}\M_u.
\end{equation}
Indeed, the linear extensions of $P$ are precisely the permutations such that
$i<_P j\Rightarrow \sigma^{-1}(i)<\sigma^{-1}(j)$.

If one specializes $A$ to the alphabet 
$A_{n+1}=\{a_0,a_1\ldots,a_n\}$,
$\Gamma_P(A_{n+1})$ becomes the sum of the integral points of $Q_P$. Their number is therefore
$E_{Q_P}(n)=\Gamma_P(n+1)$.

The change of sign $A\mapsto -A$ of the alphabet is defined on symmetric functions by means
of the $\lambda$-ring structure:
$p_n(-X)=-p_n(X)$, and one defines more generally, the multiplication of the alphabet by an element of binomial type 
 $p_n(\alpha X)=\alpha p_n(X)$.

These transformations can be naturally extended to quasi-symmetric functions. One first defines them
on $\Sym$ by setting
$\sigma_t(\alpha A)=\sigma_t(A)^\alpha$, 
then one extends to $QSym$ by defining
$\sigma_t(X\alpha \cdot A)=\sigma_t(XA)*\sigma_1(\alpha A)$. 
These transformations can then be extended to $\WQSym$ by means of the internal product of
$\WQSym^*$ \cite{NTsuper}. One obtains
\begin{equation}
\M_u(-A) = (-1)^{\max(u)}\sum_{v\le u}\M_v(A)
\end{equation}
where the sum runs over the refinement order on packed words\footnote{$v\le u$ iff the set composition
encoded by $v$ is obtained by merging adjacent blocks of that encoded by $u$.}.

If one sets $A=\{a_0,a_1,a_2\ldots\}$ et $A'=\{a_1,a_2,\ldots\}$,
one has
\begin{equation}
(-1)^n\Gamma_P(-A') = \sum_{v\in \dot{C}(P)}\M_v(A')
\end{equation}
where $\dot C(P)$ is the set of packed words satisfying
$i<_P j\Rightarrow u_i<u_j$,
otherwise said, of the packed words of the interior points of the cone.
The interior points of the polytope $nQ_P$ are obtained by evaluating on the
alphabet $\{a_1,\ldots,a_{n-1}\}$.

The number of interior points is thus
$(-1)^{n}\Gamma_P(1-n)=E_{Q_P}(-n)$, 
we have therefore in this particular case a noncommutative lift of the Ehrhart reciprocity formula.

For example,
\begin{equation}
	\begin{split}
	\Gamma_{\arbdgb}(X) &=\F_{123}+\F_{213}=\G_{123}+\G_{213}\\
		&=\M_{123}+ \M_{122}+ \M_{112}+ \M_{111}+ \M_{213}+ \M_{212} 
	\end{split}
\end{equation}
has as commutative image $F_{12}+F_3$ and as evaluation on a scalar
	$ {\alpha+2\choose 3}+{\alpha+1\choose 3}$
so that the Ehrhart polynomial of the order polytope $Q=\{0\le x_1,x_2\le x_3\}$
is
\begin{equation}
	E_Q(x)=  {x+3\choose 3}+{x+2\choose 3}  =   \frac{(x+1)(x+2)(2x+3)}{6}
\end{equation}
which is indeed the specialization $q=1$ of 
\begin{equation}
\Gamma_{\arbdgb}(X_{q,t})=
\frac{{\left(q^{3} x + q^{2} x + q^{2} + q + 1\right)} {\left(q^{2} x + q + 1\right)} {\left(q x + 1\right)}}{{\left(q^{2} + q + 1\right)} {\left(q + 1\right)}}
\end{equation}
The specialization $x=[n]_q$ gives the $q$-counting of the integral points of $nQ$ by sum of the coordinates.
Indeed, this amounts to setting $t=q^n$ in \eqref{eq:defsqx}, so that by \cite[Prop. 8.4]{NCSF2}
\begin{equation}
	\sqx \mapsto \sigma_1(X_{q,q^n}A) := \prod_{0\le i\le n}^\rightarrow\sigma_{q^i}(A) = \sum_I M_I(1,q,\ldots,q^{n})S^I,
\end{equation}
that is, 
\begin{equation}
	\Gamma_P(X_{q,q^n})=\sum_{(x_1,\ldots,x_d)\in nQ\cap\ZZ^d}q^{x_1+x_2+\cdots+x_d}.
\end{equation}
For example, the 14 integral points of $2Q$ are
$$000,\ 001,\ 011,\ 101,\ 022,\ 111,\ 012,\ 102,\ 112,\ 022,\ 202,\ 122, \ 212,\ 222$$
and $\Gamma_{\arbdgb}(X_{q,q^2})=1+q+3q^2+3q^3+3q^4+2q^5+q^6$, as expected.

Now,
\begin{equation}
	(-1)^3\Gamma_{\arbdgb}(-A)=\M_{123}+\M_{213}+\M_{112}
\end{equation}
which predicts correctly that for $n=3$ the only interior point of $3Q$ is $(1,1,2)$.
Setting  $t=q^{-n}$ in \eqref{eq:defsqx} results into
\begin{equation}
	\sqx \mapsto \sigma_1(X_{q,q^{-n}}A) := \prod_{1\le i\le n-1}^\rightarrow\lambda_{-q^{-i}}(A)
\end{equation}
so that $\Gamma_P(X_{q,q^{-n}})$ is obtained by evaluating $\Gamma_P(-A)$ on the alphabet $\{x_i=q^{-i}|i=1,\ldots,n-1\}$, 
in accordance with \cite[Theorem 2.5]{Cha4}. On our example, setting $x=[-3]_q$ yields $-q^{-4}$, corresponding to the interior point $(1,1,2)$.


\end{document}